\documentclass[12pt]{article}
\usepackage{amsmath}
\usepackage{amsfonts}
\usepackage{amssymb}

\input{tcilatex}
\begin{document}

\title{A formula for the $n^{th}$ derivative of the quotient of two functions%
}
\author{Christos Xenophontos \\
Department of Mathematics and Statistics \\
University of Cyprus\\
P.O. BOX 20537\\
Nicosia 1678 \\
Cyprus}
\maketitle

\begin{abstract}
As the title suggests, we give a formula for the $n^{th}$ derivative of a 
quotient of two functions, analogous to Leibniz's formula for the product. 
This particular note has remained unpublished since 2007 (available only my website),
however it has received several citations (e.g. \cite{eg, eg2, eg3}). As a result, it now appears
on ArXiv for future reference.
\end{abstract}

In most calculus courses one encounters Leibnitz's formula for the $n^{\text{%
th}}$ derivative of the \emph{product} of two functions, but when it comes
to the \emph{quotient} this is seldomly discussed.

Does such a formula even exist? In answering this question, we may first
attempt to answer the simpler query: what is the formula for the $n^{\text{th%
}}$ derivative of $1/f(x)$? The answer to the latter question was given in 
\cite{SV}, even though the notation used makes the formula difficult to use.
The former, and original, question was answered in \cite{G} using a linear
algebra approach which ultimately gave a rather beautiful formula involving
the determinant of a certain matrix. As the author of that article admitted,
the interest in the particular determinantal formula \textquotedblleft lies
mainly in the structured it revealed\textquotedblright\ \cite{G}. Here we
built upon the ideas in \cite{G} to obtain a more usable version of the
formula for the $n^{\text{th}}$ derivative of the quotient of two functions $%
u(x)$ and $v(x)$, that is analogous to Leibnitz's formula.

As in \cite{G}, writing $y=u/v$ and using Leibnitz's formula we have%
\begin{equation*}
u^{(r)}=\left( yv\right) ^{(r)}=\underset{j=0}{\overset{r}{\sum }}\frac{r!}{%
(r-j)!j!}y^{(j)}v^{(r-j)},
\end{equation*}%
where $u^{(r)}$ denotes $d^{r}u/dx^{r}$ for $r=0,1,2,....$ \ From the above
expression we get%
\begin{equation*}
\frac{u^{(r)}}{r!}\underset{j=0}{=\overset{r}{\sum }}\frac{y^{(j)}}{j!}\frac{%
v^{(r-j)}}{(r-j)!},
\end{equation*}%
which, for $r=0,1,2,...,n$ is nothing but the lower triangular linear system
of $n+1$ equations in the unknowns $y,\frac{y^{(1)}}{1!},\frac{y^{(2)}}{2!}%
,...,\frac{y^{(n)}}{n!},$ shown below:%
\begin{equation*}
\left[ 
\begin{array}{ccccc}
v & 0 & 0 & \ldots & 0 \\ 
\frac{v^{(1)}}{1!} & v & 0 & \ldots & 0 \\ 
\frac{v^{(2)}}{2!} & \frac{v^{(1)}}{1!} & v & \ddots & \vdots \\ 
\vdots & \vdots & \ddots & \ddots & 0 \\ 
\frac{v^{(n)}}{n!} & \frac{v^{(n-1)}}{(n-1)!} & \ldots & \frac{v^{(1)}}{1!}
& v%
\end{array}%
\right] \left[ 
\begin{array}{c}
y \\ 
\frac{y^{(1)}}{1!} \\ 
\frac{y^{(2)}}{2!} \\ 
\vdots \\ 
\frac{y^{(n)}}{n!}%
\end{array}%
\right] =\left[ 
\begin{array}{c}
u \\ 
\frac{u^{(1)}}{1!} \\ 
\frac{u^{(2)}}{2!} \\ 
\vdots \\ 
\frac{u^{(n)}}{n!}%
\end{array}%
\right] .
\end{equation*}%
In \cite{G} the above linear system was solved using Cramer's rule to obtain
a determinantal formula for $y^{(n)}$. If, instead, we take advantage of the
lower triangular structure of the linear system and use back substitution to
solve it we have%
\begin{equation*}
\frac{y^{(n)}}{n!}=\frac{1}{v}\left( \frac{u^{(n)}}{n!}-\underset{j=1}{%
\overset{n}{\sum }}\frac{v^{(n+1-j)}}{(n+1-j)!}\frac{y^{(j-1)}}{(j-1)!}%
\right) ,
\end{equation*}%
from which the following formula arises:%
\begin{equation*}
\left( \frac{u}{v}\right) ^{(n)}=\frac{1}{v}\left( u^{(n)}-n!\underset{j=1}{%
\overset{n}{\sum }}\frac{v^{(n+1-j)}}{(n+1-j)!}\frac{\left( \frac{u}{v}%
\right) ^{(j-1)}}{(j-1)!}\right) .
\end{equation*}


\begin{thebibliography}{9}

\bibitem{eg} L Barabesi, The computation of the probability density and distribution 
functions for some families of random variables by means of the Wynn-$\rho$ accelerated 
Post-Widder formula, \textit{Communications in Statistics-Simulation and Computation},
\textbf{49} (2020), 1333--1351.

\bibitem{eg2} G. Panou and R. Korakitis, The direct geodesic problem and an
approximate analytical solution in Cartesian coordinates on a triaxial ellipsoid,
\textit{Journal of Applied Geodesy}, \textbf{14} (2020), 205--213.

\bibitem{eg3} R Basu, A new formula for investigating delay integro-differential 
equations using the differential transform method involving a quotient of two functions,
\textit{Rocky Mountain Journal of Mathematics}, \textbf{51} (2021), 413--421.
\bibitem{SV} P. Shieh and K. Verghese, A general formula for the $n^{\text{th}}$ 
derivative of $1/f(x)$, \textit{Amer. Math. Monthly}, \textbf{74} (1967), 1239--1240.

\bibitem{G} F. Gerrish, A useless formula? \textit{The Mathematical Gazette},
 \textbf{64} (1980), 52.



\end{thebibliography}
\end{document}